\documentclass[centertags,12pt]{amsart}
\usepackage{latexsym}
\usepackage{amssymb}

\makeatletter
\renewcommand{\subsection}{\@startsection
{subsection}{2}{0mm}{\baselineskip}{-0.25cm}
{\normalfont\normalsize\bf}}
\makeatother

\newtheorem*{theorem*}{Theorem}

\newtheorem{proposition}{Proposition}
\newtheorem{corollary}{Corollary}
\newtheorem{lemma}{Lemma}
\theoremstyle{definition}
{\theoremstyle{remark}

\newtheorem*{claim*}{Claim}
\newtheorem{claim}{Claim}
\newtheorem{fact}{Fact}}

\def\F{\mathbb F}
\def\P{\mathbb P}
\def\N{\mathbb N}

\def\cD{\mathcal D}

\def\cH{\mathcal H}

\def\cX{\mathcal X}

\def\fq{\mathbb F_{q^2}}
\def\supp{{\rm Supp}}
\def\div{{\rm div}}
\def\dim{{\rm dim}}
\def\deg{{\rm deg}}
\def\det{{\rm det}}
\def\frx{{\rm Fr}_{\mathcal X}}

\begin{document}

\author[M.~Abd\'on]{Miriam Abd\'on}
\author[F.~Torres]{Fernando Torres}\thanks{1991 Math. Subj. Class.:
Primary 11G, Secondary 14G}

\title[Maximal curves in characteristic two]{On maximal curves in
characteristic two}
\address{IMPA, Est. Dna. Castorina 110, Rio de Janeiro, 22.460-320-RJ,
Brazil}
\email{miriam@impa.br}
\address{IMECC-UNICAMP, Cx. P. 6065, Campinas, 13083-970-SP, Brazil}
\email{ftorres@ime.unicamp.br}

\begin{abstract}

The genus $g$ of an $\fq$-maximal curve satisfies $g=g_1:=q(q-1)/2$ or
$g\le g_2:= \lfloor\text{$(q-1)^2/4$}\lfloor$. Previously,
$\fq$-maximal curves with $g=g_1$ or $g=g_2$, $q$ odd, have been
characterized up to $\fq$-isomorphism. Here it is shown that an
$\fq$-maximal curve with genus $g_2$, $q$ even, is $\fq$-isomorphic to the
nonsingular model of the plane curve $\sum_{i=1}^{t}y^{q/2^i}=x^{q+1}$,
$q=2^t$, provided that $q/2$ is a Weierstrass non-gap at some point of
the curve.

\end{abstract}

\maketitle

{\bf 1.} A projective geometrically irreducible nonsingular algebraic
curve defined over $\fq$, the finite field with $q^2$ elements, is
called {\em $\fq$-maximal} if the number
of its $\fq$-rational points  attains the Hasse-Weil upper bound
$$
q^2+1+2qg\, ,
$$
where $g$ is the genus of the curve. Maximal curves became useful in
Coding Theory after Goppa's paper \cite{goppa}, and have been intensively
studied in \cite{sti-x}, \cite{geer-vl1} (see also the references therein),
\cite{geer-vl2}, \cite{ft1}, \cite{fgt}, \cite{ft2},
\cite{gt}, \cite{chkt}, \cite{ckt1}, \cite{g-sti-x} and \cite{ckt2}.

The key property of a $\fq$-maximal curve $\cX$ is the existence of a
base-point-free linear system $\cD_{\cX}:=|(q+1)P_0|$, $P_0\in \cX(\fq)$,
defined on $\cX$ such that \cite[\S1]{fgt}
  \begin{enumerate}
\item[(1.1)] $qP+\frx(P)\in \cD_\cX$,
\item[(1.2)] $\cD_\cX$ is simple,
\item[(1.3)] $\dim(\cD_\cX)\ge 2$,
  \end{enumerate}
where $\frx$ denotes de Frobenius morphism on $\cX$ relative to $\fq$.
Then via St\"ohr-Voloch's approach to the
Hasse-Weil bound \cite{sv} one can establish arithmetical and geometrical
properties of maximal curves. In addition, Property (1.2) allows the
use of Castelnuovo's genus bound in projective spaces \cite{cas}, \cite[p.
116]{acgh}, \cite[Corollary 2.8]{ra}. In particular,
the following relation involving the genus $g$ of $\cX$ and
$n:=\dim(\cD_\cX)-1$ holds \cite[p. 34]{fgt}
    \begin{equation}\label{eq1.1}
2g\le\begin{cases}
(q-n/2)^2/n & \text{if $n$ is even}\, ,\\
((q-n/2)^2-1/4)/n & \text{otherwise}\, .
\end{cases}
    \end{equation}
It follows that
$$
g\le g_1:= q(q-1)/2\, ,
$$
which is a result pointed out by Ihara \cite{ihara}. As a matter of fact,
the so called {\em Hermitian curve}, i.e. the plane curve $\cH$ defined by
$$
Y^qZ+YZ^q=X^{q+1}\, ,
$$
is the unique $\fq$-maximal curve whose genus is $g_1$ up to
$\fq$-isomorphism \cite{r-sti}. Moreover, $\cH$ is the unique
$\fq$-maximal curve
$\cX$ such that $\dim(\cD_\cX)=2$ \cite[Thm 2.4]{ft2}.
Therefore, if $g<g_1$, then $\dim(\cD_\cX)\ge
3$ and hence (\ref{eq1.1}) implies \cite{sti-x}, \cite{ft1}
$$
g\le g_2:= \lfloor\text{$(q-1)^2/4$}\rfloor\, .
$$
If $q$ is odd, there is a unique $\fq$-maximal curve, up to
$\fq$-isomorphism, whose genus belongs to the interval
$](q-1)(q-2)/4,(q-1)^2/4]$, namely the
nonsingular model of the plane curve
$$
y^q+y=x^{(q+1)/2}\, ,
$$
whose genus is $g_2=(q-1)^2/4$ \cite[Thm. 3.1]{fgt}, 
\cite[Prop. 2.5]{ft2}.

The purpose of this paper is to extend this result to even
characteristic provided that a condition on Weierstrass non-gaps is
satisfied. 
For $q$ even, say $q=2^t$, notice that $g_2=q(q-2)/4$ and that the
nonsingular model of the plane curve
   \begin{equation}\label{eq1.2}
\sum_{i=1}^{t}y^{q/2^i} =x^{q+1}
    \end{equation}
is an $\fq$-maximal curve of genus $g_2$.

     \begin{theorem*}
Let $q$ be even, $\cX$ a $\fq$-maximal curve of genus $g$ having
both properties:
\begin{enumerate}
\item $(q-1)(q-2)/4<g\le g_2=q(q-2)/4$, and
\item There exists $P\in \cX$ such that $q/2$ is a Weierstrass non-gap at
$P$. 
\end{enumerate}
Then $\cX$ is $\fq$-isomorphic to the nonsingular model of the plane curve
defined by Eq. (\ref{eq1.2}). In particular, $g=g_2$.
     \end{theorem*}

Let $\cX$ be a $\fq$-maximal of genus $g\in ](q-1)(q-2)/4, q(q-2)/4]$,
$q$ even, and $P\in \cX$. We have that $P\in\cX(\fq)$ if $q/2$
is a Weierstrass non-gap at $P$, see Corollary \ref{cor2.2}. Now, on the
one hand,  from Corollary \ref{cor2.1}, $\cX$ only admits two types of
Weierstrass semigroups at $\fq$-rational points, namely either semigroups
of type $\langle q/2,q+1\rangle$ or semigroups of type $\langle
q-1,q,q+1\rangle $. On the other hand, from Proposition \ref{prop5.1},
$\cX$ satisfies the second hypothesis of the theorem provided that
$\cX$ is $\fq$-covered by $\cH$. Therefore, if there existed a $\fq$-maximal
curve of genus $g\in ](q-1)(q-2)/4,q(q-2)/4]$ for which the Weierstrass
semigroup at any $\fq$-rational point is $\langle q-1,q,q+1\rangle $,
then such a curve could not
be $\fq$-covered by the Hermitian curve. As far as we know, the
existence of
maximal curves not covered by the Hermitian is an open problem. We notice
that the nonsingular model $\cX$ of
the plane curves
$y^q+y=x^m$, $m$ a divisor of $q+1$, have been characterized as those
curves such that $m_1(P)n=q+1$ for some $P\in \cX(\fq)$, where $m_1(P)$
stands for the first positive Weierstrass non-gap at $P$ and
$n=\dim(\cX)-1$; see \cite[\S2]{fgt}. Moreover, the hypothesis on
Weierstrass non-gaps cannot be relaxed, cf. \cite[p. 37]{ft2}, 
\cite[Remark 4.1(ii)]{chkt}.

We prove the theorem by using some properties of maximal curves stated in
\cite{fgt}, \cite{ft2} and \cite{ckt1}, Castelnuovo's genus bound in
projective spaces, and Frobenius orders which were introduced by
St\"ohr and Voloch \cite{sv}. For basic facts on Weierstrass point theory
and Frobenius orders the reader is referred to \cite{sv}.


{\bf 2. Proof of the Theorem.} Let $\cX$ be a $\fq$-maximal curve of genus
$g$ large enough. The starting point of the proof is the computation of
some invariants for the following linear systems:
$$
\cD:=\cD_{\cX}=|(q+1)P_0|\quad \text{and}\quad 2\cD:=
2\cD_{\cX}=|2(q+1)P_0|\, ,
$$
where $P_0\in \cX(\fq)$. For $P\in \cX$, let $(m_i(P): i\in \N_{0})$
be the strictly increasing sequence that enumerates the Weierstrass
semigroup $H(P)$ at $P$.

   \begin{lemma}\label{lemma2.1} Let $\cX$ be a $\fq$-maximal curve of
genus $g$ such that
$$
(q-1)(q-2)/4 < g\le (q-1)^2/4\, .
$$
Then the following properties hold:
   \begin{enumerate}
\item We have $\dim(\cD)=3$.
\item If $P\in \cX(\fq)$, then the $(\cD,P)$-orders are $0,1, q+1-m_1(P)$
and $q+1$. If $P\not\in \cX(\fq)$, then the set of $(\cD,P)$-orders
contains the elements $q-m_i(P)$, $i=0,1,2$.
\item We have $\dim(2\cD)=8$.
   \end{enumerate}
   \end{lemma}
   \begin{proof} (1) From Iq. (\ref{eq1.1}) and the lower bound on
$g$ it follows that $\dim(\cD)\le 3$ (indeed, we obtain this result for
$(q-1)(q-2)/6<g $. If we had $\dim(\cD)=2$, then from \cite[Thm. 2.4]{ft2}
it would follow $g=q(q-1)/2$, contradiction. Thus $\dim(\cD)=3$.

(2) See \cite[Prop. 1.5(ii)(iii)]{fgt}.

(3) An easy computation shows that $2m_3(P_0)\ge 8$, since $m_2(P_0)=q$
and $m_3(P_0)=q+1$ \cite[Prop. 1.5(iv)]{fgt}. Hence $\dim(2 \cD)\ge 8$;
the equality follows  from Castelnuovo's genus bound and the lower bound
on $g$.
\end{proof}

\begin{corollary}\label{cor2.1} Let $\cX$ be as in the previous lemma
and suppose that $q$ is even, $q>4$.

  \begin{enumerate}
\item For $P\in \cX(\fq)$,
\subitem(i) the $(\cD,P)$ orders are either $0,1,2,q+1$ or $0,1,
q/2+1,q+1$;
\subitem(ii) either $m_1(P)=q-1$ or $m_1(P)=q/2$. Equivalently, the first
three positive Weierstrass non-gaps at $P$ are either $q-1, q,q+1$  or $q/2,q,
q+1$.
\item For $P\not\in \cX(\fq)$, the $(\cD,P)$-orders are either
$0,1,2,q$ or $0,1,q/2,q$.
   \end{enumerate}
\end{corollary}

   \begin{proof} (1) Let $P\in \cX(\fq)$ and set $j:=q+1-m_1(P)$. By Lemma
\ref{lemma2.1}(2), it is enough
to prove (i). From that result and the definition of
$2\cD$, the following set
$$
O:=\{ 0, 1, 2, j, j+1, 2j, q+1, q+2, q+1+j, 2q+2\}
$$
is contained in the set of $(2\cD,P)$-orders. Since $\dim(2\cD)=8$ (see
Lemma \ref{lemma2.1}(3)), $\# O \le 9$. We observe that $j<q$: otherwise
$g=0$. So if $j>2$, then $2j=q+2$, as $q$ is even, and the result
follows.

(2) From Lemma \ref{lemma2.1}(1) and \cite[Thm. 1.4(ii)(iii)]{fgt}, the
$(\cD,P)$-orders are $0,1,j=j(P)$ and $q$ with $2\le j\le q-1$. We claim
that $j<q-1$. Otherwise $(q-1)P+D_P\sim qP+\frx(P)$, with
$P\not\in\supp(D_P)$ and so $\cX$ would be hyperelliptic; then
$1+q^2+2qg\le 2(1+q^2)$ and hence $(q-1)(q-2)/4<g\le q/2$, a contradiction
since we have assumed $q>4$. Now, the following set
$$
\{0,1,2,j,j+1,2j,q,q+1,q+j,2q\}
$$
is contained in the set of $(2\cD,P)$-orders and the result follows as in
the proof of item (1).
   \end{proof}
 \begin{corollary}\label{cor2.2} Let $\cX$ be as in Lemma \ref{lemma2.1}
and suppose that $q$ is even, $q>4$. Suppose also that $q/2$ is a
Weierstrass non-gap at $P\in \cX$. Then $P\in \cX(\fq)$.
 \end{corollary}
  \begin{proof} Suppose that $P\not\in\cX(\fq)$. Let $x\in \bar\fq(\cX)$
such that $\div_\infty(x)=q/2P$. Let $e:=v_{\frx(P)}(x-x(\frx(P))$. Then
$\div(x-x(\frx(P))=e\frx(P)+D-q/2P$ with $P, \frx(P)\not\in \supp(D)$.
>From Property (1.1), both $e+1$ and $2e+1$ are $(\cD,\frx(P))$-orders.
Since $\frx(P)\not\in\cX(\fq)$ and $e\ge 1$, from Corollary
\ref{cor2.1}(2) follows that $q=3$, a contradiction.
  \end{proof}

Now, for $q>4$ the theorem follows from the proposition below. The case
$q=4$ is considered in \S4.

    \begin{proposition}\label{prop2.1} Let $\cX$ be a projective
geometrically irreducible nonsingular algebraic curve over $\fq$, $q$
even. The  following statements are equivalent:
     \begin{enumerate}
\item The curve $\cX$ is $\fq$-isomorphic to the non-singular model of the
plane curve given by Eq. (\ref{eq1.2}).
\item The curve $\cX$ is $\fq$-maximal of positive genus,
$\dim(\cD_{\cX})=3$, and there
exists $P_0\in \cX(\fq)$ such that $q/2$ is a Weierstrass non-gap at
$P_0$.
\item The curve is $\fq$-maximal and there exists $P_1\in \cX(\fq)$ such
that for $\cD:=\cD_{\cX}=|(q+1)P_1|$, the following holds:
\subitem(i) the $(\cD,P_1)$-orders are $0,1, q/2+1, q+1$;
\subitem(ii) the $(\cD,P)$-orders are $0,1,2,q+1$ if $P\in
\cX(\fq)\setminus \{P_1\}$;
\subitem(iii) the $(\cD,P)$-orders are $0,1,2,q$ if $ P\in \cX\setminus
\cX(\fq)$;
\subitem(iv) the $\fq$-Frobenius orders of $\cD$ are $0,1,q$.
    \end{enumerate}
    \end{proposition}

%
%

{\bf 3. Proof of Proposition \ref{prop2.1}.} Throughout this section we
assume $q\ge 4$ since the case $q=2$ is trivial.

$(1)\Rightarrow (2):$ The non-singular
model of (\ref{eq1.2}) is $\fq$-covered by
the Hermitian curve, and so it is maximal by \cite[Prop. 6]{lachaud}. The
unique point $P_0$ over $x=\infty$ is $\fq$-rational and $q/2$ and $q+1$
are Weierstrass non-gaps at $P_0$. Since the genus of the curve is
$q(q-2)/2$, it follows that $\dim(|(q+1)P_0|)=3$.

$(2)\Rightarrow (3):$ This implication is a particular case of
\cite[p. 38]{ft2}; for the sake of completeness we write the proof.  Take
$P_1=P_0$. Then $m_1(P_0)=q/2$,
$m_2(P_0)=q$ and $m_3(P_0)=q+1$, cf. Corollary \ref{cor2.1}(1)(ii).
The case $P=P_0$ follows from Lemma \ref{lemma2.1}(2). Let $P\in
\cX\setminus\{P_0\}$. Let $x\in \fq(\cX)$ such that
$\div_\infty(x)=m_1(P_0)P_0$. Then $e_P:=v_P(x-x(P))$ and $2e_P$ are
$(\cD,P)$-orders. We claim that $e_P=1$; otherwise $0, 1, e_P$, and $2e_P$
would be $(\cD,P)$-orders and hence, by \cite[Thm. 1.4(ii)]{fgt} and
being $q$ even, we would have $e_P=q/2$ and $P\not\in\cX(\fq)$. Therefore
$q/2P\sim q/2P_0$ and from Property (1.1) (and since the genus of $\cX$
is positive) we would have $\frx(P)=P_0$, a contradiction. Thus, by
\cite[Thm. 1.4(ii)(iii)]{fgt}, the $(\cD,P)$-orders are $0,1,2$
and $q+1$ (resp. $0,1,2,$ and $q$) if $P\in \cX(\fq)$ (resp. $P\not\in
\cX(\fq)$). Finally, the
assertion on $\fq$-Frobenius orders follows from $\dim(\cD)=3$ and
\cite[\S2.2]{ft2}.

$(3)\Rightarrow (1)$: By Lemma \ref{lemma2.1}(2), $m_1(P_1)=q/2$,
$m_2(P_1)=q$ and $m_3(P_1)=q+1$. Let $x, y\in \fq(\cX)$ such that
$$
\div_\infty(x)=q/2P_1\, ,\qquad\text{and}\qquad\div_\infty(y)=(q+1)P_1\, .
$$
Then $\cX$ admits a $\fq$-plane model of type
   \begin{equation}\label{eq4.1}
x^{q+1} + a y^{q/2} + \sum_{i=0}^{q/2 -1} A_i(x) y^i=0\, ,
   \end{equation}
where $a\in \fq^*$ and  $A_i(x) \in \fq[x]$ with $\deg(A_i(x)) \le q-2i$,
$i=0, \dots, q/2 -1$. This equation is usually referred to as the {\em
Weierstrass canonical form over $\fq$ of $\cX$}, see e.g. \cite[Lemma
3]{k} and the references therein.

Next we use $x$ as a separating variable of $\fq(\cX) \!\mid \!
\fq$, and denote by $D^i:=D^i_x$ the $i$th Hasse derivative with respect
to $x$. Properties of
these operators can be found e.g. in \cite[\S3]{hefez}. In particular, we
recall the following facts: For $z, w \in \bar\fq(\cX)$,
  \begin{enumerate}
\item[(H1)]\quad $D^i(z+w)=D^i(z)+D^i(w)$,
\item[(H2)]\quad $D^i(zw)=\sum_{j=0}^{i}D^{i-j}(z)D^j(w)$,
\item[(H3)]\quad
$$
D^iz^{2j}=\begin{cases}
(D^{i/2}z^j)^2 & \text{if $i$ is even}\, ,\\
0         & \text{otherwise}\, .
\end{cases}
$$
   \end{enumerate}
Then, for $q'$ a power of two, (H3) implies:
   \begin{enumerate}
\item[(H3')]
$$
D^iz^{q'} = \begin{cases}
(D^{i/q'} z)^{q'} & \text{if $i\equiv 0\pmod{q'}$}\, ,\\
0                & \text{otherwise}\, .
\end{cases}
$$
   \end{enumerate}
Now, the morphism associated to $\cD$ is given by $(1:x:x^2:y)$. Since
the $\cD$-orders are $0,1,2$ and $q$, for $i=3,\ldots,
q-1$, we have (see \cite[p. 5]{sv})
\begin{align*}
\det
\begin{pmatrix}1 & x & x^2 & y \\
 0 & 1 & 0 & Dy \\
 0 & 0 & 1 & D^2y \\
 0 & 0 & 0 & D^iy
\end{pmatrix}  & = D^iy = 0\, .\\
\intertext{We also have}
\det
\begin{pmatrix}1 & x^{q^2} & x^{2q^2} & y^{q^2} \\
 1 & x & x^2 & y \\
 0 & 1 & 0 & Dy \\
 0 & 0 & 1 & D^2y \\
\end{pmatrix} & = 0\, ,
\end{align*}
or equivalently,
     \begin{equation}\label{eq4.2}
y+y^{q^2}+(x+x^{q^2})Dy+(x^2+x^{2q^2})D^2y=0\, ,
     \end{equation}
since the $\fq$-Frobenius orders of $\cD$ are $0, 1$ and $q$ (see
\cite[Prop. 2.1]{sv}).

\begin{claim}\label{claim3.1} Eq. (\ref{eq4.1}) can be simplified to
\begin{equation}\label{eq4.21}
\sum_{i=1}^{t}a_i y^{q/2^i}+b=x^{q+1}\, ,
\end{equation}
where $a_1, \ldots , a_t, b \in \bar \fq$ with $a_t\in \fq$.
\end{claim}

Let us first show how this claim implies Proposition \ref{prop2.1}(1). To
do so, let  $\alpha \in \bar\F_{q^2}$ such that
$$
\sum_{i=1}^{t} a_i\alpha^{q/2^i} = b \, .
$$
Then, with  $z:= y + \alpha$, the curve  $\cX$ is
$\bar\F_{q^2}$-isomorphic to the non-singular model of the curve defined
by
$$
\sum_{i=1}^{t} a_i z^{q/2^i}=x^{q+1}\, .
$$

\begin{fact}\label{fact0} The element $a_t$ can be assumed to be equal to
one. If so, then
\begin{enumerate}
\item $a_{t-1}=a_1^{-2}$.
\item $a_i=a_{t-1}^{2^{t-i}-1}$, $i=1,\ldots , t-1$. In particular, $a_i
\in \fq^*$ for each $i$.
\item $\alpha \in \fq$.
\end{enumerate}
\end{fact}

\begin{proof}({\em Fact \ref{fact0}}) From Eq. (\ref{eq4.21}), we have
$$
a_tDy=x^q \qquad \text{and} \qquad a_t^3 D^2y=a_{t-1} x^{2q}\, ,
$$
and so $a_t \neq 0$. Hence we can assume $a_t=1$ via the automorphism
$(x,y)\mapsto (x, a_t y)$. Now from Eq. (\ref{eq4.2}) we obtain
$$
y+y^{q^2}+ x^{q+1}+ (x^{q+1})^q + a_{t-1}( (x^{q+1} )^2+
(x^{q+1})^{2q})=0\, .
$$
This relation together with  Eq. (\ref{eq4.21}) imply
\begin{align*}
{} & (1 + a_{t-1} a_{1}^{2q}) y^{q^2} + \sum_{i=1}^{t-1} {(a_i^q+a_{t-1}
a_{i+1}^{2q}) y^{{q^2}/2^{i}}} +(1+a_{t-1}a_{1}^2 ) y^q +\\
{} & \sum_{i=1}^{t-1} {a_i+(a_{t-1} a_{i+1}^{2})
y^{{q}/2^{i}}}+b+b^q+a_{t-1}(b^2+b^{2q})=0\, .
\end{align*}
Therefore, as $v_{P_1}(y)<0$, the following identities hold
   \begin{enumerate}
\item[(i)]\quad $1+a_{t-1} a_1^{2q}=0$,
\item[(ii)]\quad $1+a_{t-1} a_1^2=0$,
\item[(iii)]\quad $a_i+a_{t-1} a_{i+1}^{2}= 0$,\,  $i=1,\ldots , t-1$,
\item[(iv)]\quad $a_i^q+a_{t-1} a_{i+1}^{2q}= 0$,\, $ i=1, \ldots,
t-1$,\quad and,
\item[(v)]\quad $b+b^q+ a_{t-1}(b^2 +b^{2q})=0$.
   \end{enumerate}

>From (i) and (iii) follow Items 1 and 2. To see Item 3 we replace
$b by \sum_{i=1}^{t}a_i\alpha^{q/2^i}$ in (v). After some computations and
using (i)--(iv) we find that $\alpha+\alpha^{q^2}=0$ and
the proof of Fact \ref{fact0} is complete.
\end{proof}

Consequently, the automorphism $(x,y)\mapsto (x,y+\alpha)$ is indeed
defined over $\fq$. Finally let $x_1:= a_1^{-1}x$ and $y_1:=a_{t-1}z$.
Then from Fact \ref{fact0} we obtain
$$
\sum_{i=1}^{t} y_1^{q/2^i}=x_1^{q+1}\, ,
$$
which shows Proposition \ref{prop2.1}(1).

({\em Proof of Claim \ref{claim3.1}.}) Suppose that $x$ and $y$ satisfy a
relation of type
    \begin{equation}\label{eq4.3}
x^{q+1}+ay^{q/2}+\sum_{i=0}^{2^{s-2}-1} A_{i2^{t+1-s}}(x)y^{i2^{ t+1-s}}+
\sum_{i=s}^{t}a_i y^{2^{t-i}}=0\, ,
    \end{equation}
where $a\in \fq^*$, $2\le s\le t+1$, and $a_i\in \fq$ for each $i$. Recall
that $q=2^t$ and notice that Eq. (\ref{eq4.1}) provides such a relation
for $s=t+1$.
     \begin{fact}\label{fact1} For $2\le s\le t+1$, we have
\begin{enumerate}
\item \quad $A_{(2i+1)2^{t+1-s}}(x)=0$,\ \  $i\ge 1$.
\item \quad $A_{2^{t+1-s}}(x)\in \fq$.
\end{enumerate}
    \end{fact}
\begin{proof} ({\em Fact \ref{fact1}}) By applying $D^{2^{t+1-s}}$ to Eq
(\ref{eq4.3}) and using properties (H1)--(H3) and (H3') above we have that
\begin{equation}\label{eq4.4}
a'\Gamma+\sum_{i=0}^{2^{s-2}-1}\Delta_i +
\sum_{i=0}^{2^{s-2}-1}A_{i2^{t+1-s}}(x)(Dy^i)^{2^{t+1-s}}+
\Lambda+\Psi=0\, ,
\end{equation}
where
  \begin{align*}
\Gamma:= & \begin{cases}
(Dy)^{2^{t-1}} & \text{if $s=2$}\, ,\\
0              & \text{otherwise}\, ;
\end{cases}  \\
\Delta_i:= & \begin{cases}
D(A_i(x))y^i  & \text{if $s=t+1$}\, ,\\
\sum_{j=0}^{2^{t-s}-1}D^{2^{t+1-s}-2j}(A_{i2^{t+1-s}}(x))D^{2j}(y^{i2^{t+1-s}})
&  \text{otherwise}\, ;
\end{cases} \\
\Lambda:= & \begin{cases}
0   &  \text{if $s=t+1$}\, ,\\
(D^2y)^{2^{t-s}} &  \text{otherwise}\, ;
\end{cases} \\
\intertext{and}
\Psi:=\begin{cases}
x^q  & \text{if $s=t+1$}\, ,\\
0    & \text{otherwise}\, .
\end{cases}
  \end{align*}

Since $D^{2^{t+1-s}}(y^{i2^{t+1-s}})=(Dy^i)^{2^{t+1-s}}$, Eq.
(\ref{eq4.4}) becomes
\begin{equation}\label{eq4.5}
a'\Gamma+\sum_{i=0}^{2^{s-2}-1}\Delta_i + F (Dy)^{2^{t+1-s}} +
\Lambda+\Psi=0\, ,
\end{equation}
where
$$
F:= \sum_{i=0}^{2^{s-2}-1}
A_{(2i+1)2^{t+1-s}}(x)y^{i2^{t+2-s}}\, .
$$
Next we show that $v_{P_1}(F)=0$ $(*)$. This will imply Fact \ref{fact1}
since
$$
v_{P_1}(F)=\min\{ v_{P_1}(A_{(2i+1)2^{t+1-s}}(x)y^{i2^{t+2-s}}): i=0,
\ldots,
2^{s-2}-1\}\, .
$$
To see $(*)$, we first compute $v_{P_1}(Dy)$ and $v_{P_1}(D^2y)$.
For a local parameter $t$ at $P_1$, we have
$$
v_{P_1}(Dy) = v_{P_1}(dy/dt) - v_{P_1}(dx/dt)=-q-2-v_{P_1}(dx/dt)\, .
$$
To calculate $v_{P_1}(dx/dt)$, we use the fact that $x:\cX\to
\P^1(\bar\fq)$ is totally ramified at $P_1$, and that $v_P(x-x(P))=1$ for
each  $P\in \cX\setminus\{P_1\}$ (cf. proof of $(2)\Rightarrow (3)$). We have
then $v_{P_1}(dx/dt)=2g-2=q^2/2-q-2$ and so
$$
v_{P_1}(Dy)=-q^2/2\, .
$$
Now we compute  $v_{P_1}(D^2y)$ from Eq. (\ref{eq4.2}). In fact, as
$$
v_{P_1}(y +  y^{q^2}) = -q^2(q+1) < v_{P_1}(Dy(x + x^{q^2})) =
-q^2(q+1)/2\, ,
$$
then $v_{P_1}((x^2 + x^{2q^2})D^2y)= v_{P_1}(y + y^{q^2})$ and so
$$
v_{P_1}(D^2y) = -q^2\, .
$$
If $s=t+1$, then Eq. (\ref{eq4.5}) reads
$$
\sum_{i=0}^{q/2-1} D(A_i(x))y^i+F Dy+x^q=0\, .
$$
Thus we have
$$
v_{P_1}(F)+ v_{P_1}(Dy)=v_{P_1}(x^q)\, ,
$$
since $v_{P_1}(\sum_{i=0}^{q/2-1} D(A_i(x))y^i)> v_{P_1}(x^q)$ as
$\deg(A_i(x))\le q-2i$ for each $i$. Then $(*)$ follows because
$v_{P_1}(x^q)=-q^2/2$.

Let $3\le s\le t$. Then Eq. (\ref{eq4.5}) reads
$$
\sum_{i=0}^{ 2^{s-2}-1} \Delta_i +F(Dy)^{2^{t+1-s}}+(D^2y)^{2^{t-s}}=0\, .
$$
Consequently, as
$v_{P_1}(\sum_{i=0}^{2^{s-2}-1}\Delta_i)>2^{t-s}v_{P_1}(D^2 y)$, we have
that
$$
v_{P_1}(F)+ 2^{t+1-s}v_{P_1}(Dy)= 2^{t-s}v_{P_1}(D^2y)\, ,
$$
and the proof follows for $s\ge 3$.

Finally, let $s=2$. Then Eq. (\ref{eq4.5}) reads
$$
a'(Dy)^{q/2}+ \Delta_0 + F(Dy)^{q/2}+ (D^2y)^{q/4}=0\, .
$$
Then as above we have
$$
v_{P_1}(F) + qv_{P_1}(Dy)/2= qv_{P_1}(D^2y)/4\, ,
$$
and the proof of Fact \ref{fact1} is complete.
\end{proof}

Applying Fact \ref{fact1} for $s=t+1,\ldots, 2$, we reduce Eq.
\ref{eq4.1} to

    \begin{equation}\label{eq4.6}
x^{q+1}+ay^{q/2}+A_0(x)+\sum_{i=2}^{t}a_iy^{q/2^i}=0\, ,
     \end{equation}

where $A_0(x)=\sum_{i=0}^{q}b_ix^i\in \fq[x]$. Moreover, we can assume
$a_t=1$.

     \begin{fact}\label{fact2}
   \begin{enumerate}
\item $A_0(x)=b_0+\sum_{i=1}^{t} b_i x^{q/2^i}$.
\item $b_i=a_ib_t^{q/2^i}$, $i=1,\ldots, t$.
   \end{enumerate}
     \end{fact}
\begin{proof} ({\em Fact \ref{fact2}}) (1) Via the $\fq$-map $x\mapsto x+b_q$
applied to Eq. (\ref{eq4.6}), we
can assume $b_q=0$. Let $i$ be a natural number which is not a
power of two and satisfies the condition $3\le i<q$. Applying $D^i$
to Eq. (\ref{eq4.6}) we have (1) as $D^iy=0$.

(2) From Item (1) and Eq. (\ref{eq4.6}) we have the following equation:
   \begin{equation}\label{eq4.7}
x^{q+1}+\sum_{i=1}^{t}a_iy^{q/2^i}+
\sum_{i=1}^{t}b_ix^{q/2^i}+b_0=0\, ,
   \end{equation}
where $a_1:=a\neq 0$. Then, $Dy=x^q+b_t$ and $D^2y=a_{t-1}x^{2q}+
a_{t-1}b_t^2+b_{t-1}$. Now we want to use Eq.
(\ref{eq4.2}); thus we first have to compute $y^{q^2}+y$. Eq.
(\ref{eq4.7}) allows us to do the following computations
\begin{align*}
a_1^{2q}y^{q^2} & =
a_2^{2q}y^{q^2/2}+\sum_{i=2}^{t-1}a_{i+1}^{2q}y^{q^2/2^{i}}+x^{2q(q+1)}+
\sum_{i=0}^{t-1}b_{i+1}^{2q}x^{q^2/2^{i}}+b_0^{2q}\, ,\\
\intertext{and}
a_1^{q}y^{q^2/2}& = \sum_{i=2}^{t}a_i^{q}y^{q^2/2^{i}}+ x^{q(q+1)}+
\sum_{i=1}^{t}b_i^{q}x^{q^2/2^{i}}+b_0^{q}\, ,
\end{align*}
so that
    \begin{align*}
a_1^{3q}y^{q^2}& =
\sum_{i=2}^{t-1}(a_2^{2q}a_i^q+a_1^qa_{i+1}^{2q})y^{q^2/2^i}+
a_2^{2q}y^q+a_1^qx^{2q(q+1)}+ a_2^{2q}x^{q(q+1)}+
a_1^qb_1^{2q}x^{q^2}+\\
{} &\ \ \ \
\sum_{i=1}^{t-1}(a_2^{2q}b_i^q+a_1^qb_{i+1}^{2q})x^{q^2/2^i}+a_2^{2q}b_t^qx^q+
a_2^{2q}b_0^q+a_1^qb_0^{2q}\, .
    \end{align*}
Now, applying $D^{q/2^i}$, $1\le i<t$, to Eq. (\ref{eq4.7}) and taking
into account that $D^\ell=0$ for $3\le \ell<q$, we have
$$
a_i(Dy)^{q/2^i}+a_{i+1}(D^2y)^{q/2^{i+1}}=b_i\, ,
$$
and so, for $1\le i<t$,
\begin{enumerate}
\item[($*1$)] $a_i=a_{t-1}^{2^{t-i}-1}$, and
\item[($*2$)] $b_i=a_{i+1}b_{t-1}^{q/2^{i+1}}$.
\end{enumerate}
It follows that
$$
a_2^{2q}a_i^q+a_1^qa_{i+1}^{2q}=
a_2^{2q}b_i^q+a_1^qb_{i+1}^{2q}=0\, \qquad (i=1, \ldots, t-1)\, ,
$$
and hence
$$
a_1^{3q}y^{q^2} =
a_2^{2q}y^q+a_1^qx^{2q(q+1)}+ a_2^{2q}x^{q(q+1)}+
a_1^qb_1^{2q}x^{q^2}+
a_2^{2q}b_t^qx^q+
a_2^{2q}b_0^q+a_1^qb_0^{2q}\, .
$$
Now, from Eq. (\ref{eq4.7}) together with the following identities
coming from $(*1)$ and $(*2)$,
   \begin{itemize}
\item\quad $a_1^{3q+2}=a_2^{2q}$,
\item\quad $a_i=a_{i+1}^2a_{t-1}$, $i=1,\ldots, t-1$,
\item\quad $b_i=b_{i+1}^2a_{t-1}$, $i=1,\ldots,t-2$,
   \end{itemize}
we deduce the identity between polynomials in $x$:
\begin{equation*}
\begin{split}
 & a_1^{q+2}a_{t-1}x^{2q(q+1)}
+a_1^2a_2^{2q}a_{t-1}x^{q(q+1)}+
a_1^{q+2}a_{t-1}b_1^{2q}x^{q^2}+a_2^{2q}a_{t-1}x^{2q+2}+
 a_2^{2q}x^{q+1}+\\
& a_{t-1}(a_2^{2q}b_1^2+a_1^2a_2^{2q}b_t^q)x^q +
a_2^{2q}(a_{t-1}b_t^2+b_{t-1})x^2+
a_2^{2q}b_tx+
a_2^{2q}(a_{t-1}b_0^2+b_0)+a_1^2a_2^{2q}a_{t-1}b_0^q+\\
& a_1^{q+2}a_{t-1}b_0^{2q}
=
a_2^{2q}(x^{q^2}+x)(x^q+b_t)+a_2^{2q}(x^{2q^2}+x^2)(a_{t-1}x^{2q}+
a_{t-1}b_t^2+b_{t-1})\, .
\end{split}
\end{equation*}
Then, from the coefficients of $x^{q^2}$ we obtain
$a_1^{q+2}b_1^{2q}a_{t-1}=a_2^{2q}b_t$. Moreover, since $a_2\neq 0$ and
$a_{t-1}\in \fq$, from
$(*1)$ and $(*2)$ we have that $b_{t-1}=a_{t-1}b_t^2$. In addition, from
$(*2)$ we also have that $b_i=b_{t-1}^{q/2^{i+1}}a_{t-1}^{q/2^{i+1}-1}$,
$i=1,\ldots, t-1$, and Item (2) follows.
\end{proof}
To finish the proof of Claim \ref{claim3.1}, we apply the $\fq$-map
$(x,y)\mapsto (x,b_tx+y)$ to Eq. (\ref{eq4.7}). We obtain
a relation of type
$$
x^{q+1}+\sum_{i=1}^{t}a_iy^{q/2^i}+b_0+
\sum_{i=1}^{t}(a_ib_t^{q/2^i}+b_i)x^{q/2^i}=0\, .
$$
and Claim \ref{claim3.1} follows Fact \ref{fact2}(2).
\smallskip

%
%

{\bf 4. Case $q=4$.} Here, for $q=4$, we prove the
theorem without the hypothesis on Weierstrass non-gaps.

\begin{proposition}\label{prop4.1}
An $\F_{16}$-maximal curve $\cX$ of genus $g=2$ is $\F_{16}$-isomorphic to
the non-singular model of $y^2+y=x^5$. In particular,
it is $\F_{16}$-covered by the Hermitian curve over $\F_{16}$.
   \end{proposition}

   \begin{proof} We show that $\cX$ satisfies the hypothesis in
Proposition \ref{prop2.1}(2). Clearly, $\dim(\cD_\cX)=3$ and
$m_1(P)\in \{2,3\}$ for $P\in \cX(\F_{16})$. Suppose that
$m_1(P)=3$ for each $P\in \cX(\F_{16})$. Then by Lemma
\ref{lemma2.1}(2), the $(\cD,P)$-orders (resp. $\cD$-orders) are
0,1,2 and 5 (resp. 0,1,2 and 4). In addition, for $Q\not\in\cX(\F_{16})$,
the $(\cD,Q)$-orders are either 0,1,2 and 4 or 1,2,3 and 4. Then,
the following statements hold:
  \begin{enumerate}
\item the $2\cD$-orders are 0,1,2,3,4,5,6,7, and 8;
\item for $P\in \cX(\F_{16})$, the $(2\cD,P)$-orders are
0,1,2,3,4,5,6,7 and 10;
\item for $P\in \cX\setminus \cX(\F_{16})$, the
$(2\cD,P)$-orders are 0,1,2,3,4,5,6,7, and 8.
   \end{enumerate}
Thus $\supp(R) = \cX(\fq)$ and $v_P(R)=2$ for each $P \in
\cX(\fq)$, with $R$ being the ramification divisor associated to $2\cD$.
Thus
$$
36(2g-2)+40=\deg(R)=2\#\cX(\F_{16})=2(4(2g-2)+25)\, ,
$$
which implies $28(2g-2)=10$, a contradiction.
    \end{proof}


{\bf 5.} An $\fq$-maximal curve $\cX$, which  satisfies the hypothesis of the
theorem,
is $\fq$-isomorphic to $ \cH / \langle \tau\rangle$, where $\cH$ is the
Hermitian curve and $\tau$ an involution on $\cH$. Conversely, let us consider a
separable $\fq$-covering of curves
$$
\pi: \cH\to \cX\, .
$$
Notice that $\cX$ is $\fq$-maximal by \cite[Prop. 6]{lachaud}. Let $g$ be
the genus of $\cX$. We have the following

    \begin{proposition}\label{prop5.1} In the above situation, suppose
that $g>\lfloor\mbox{$\frac{q^2-q+4}{3}$}\rfloor$. Then
$\deg(\pi)=2$ and $g=\lfloor\mbox{$\frac{(q-1)^2}{4}$}\rfloor$. In
addition,
     \begin{enumerate}
\item $\cX$ is the non-singular model of
$y^q+y=x^{(q+1)/2}$ provided that $q$ odd
\item $\cX$ is the non-singular model of $\sum_{i=1}^{t}y^{q/2^i}=x^{q+1}$
provided that $q=2^t$.
     \end{enumerate}
     \end{proposition}
     \begin{proof} Th claim $\deg(\pi)=2$ follows from the Riemann-Hurwitz
formula taking into account the hypothesis on $g$. Then $\pi$ has (totally) ramified
points: this follows from the Riemann-Hurwitz formula and \cite{ft1}.
On the other hand, the hypothesis on $g$ allows us to use \cite[Lemma
3.1]{ckt1} and \cite[Prop. 1.5]{fgt} to conclude that
$$
m_1(P)<m_2(P)\le q<m_3(P)\qquad\text{for each}\ P\in \cX\, .
$$
Let $Q_0\in \cH$ be totally ramified for
$\pi$ and set $P_0:= \pi(Q_0)$. Then the Weierstrass non-gaps at $Q_0$
less than or equal to $2q$ are
$$
\text{either}\qquad q, 2q-1,2q\qquad\text{or}\qquad q,q+1,2q\, .
$$
It follows that $m_2(P_0)=q$ and that
$2m_1(P_0)\in\{q, q+1\}$ (see e.g. \cite[proof of Lemma 3.4]{t}). Now if
$q$ is odd,
$m_1(P_0)=(q+1)/2$ so that $m_3(P_0)=q+1$; hence $P_0\in \cX(\fq)$ and
(1) follows from \cite[Thm. 3.1]{fgt}.

If $q$ is even, we claim that $P_0$ can be chosen in $\cX(\fq)$. To see
this, as $\deg(\pi)=2$, for each $Q\in \cX(\fq)$ the product formula
gives the following possibilities:
    \begin{enumerate}
\item $\#\pi^{-1}(Q)=2$ and $\pi^{-1}(Q)\subseteq \cH(\fq)$;
\item $\#\pi^{-1}(Q)=1$ and $\pi^{-1}(Q)\in \cH(\fq)$;
\item $\#\pi^{-1}(Q)=1$ and $\pi^{-1}(Q)\in \cH(\F_{q^4})$.
    \end{enumerate}
Since $\cH(\F_{q^4})=\cH(\fq)$ we have that (3) is empty if $Q\in
\cH(\F_{q^4})\setminus \cH(\fq)$. Since $2\#\cX(\fq)>\#\cH(\fq)$, the
claim follows.

Finally we get $2m_1(P)=q$ and the result follows from Proposition
\ref{prop2.1}.
     \end{proof}
\smallskip

{\bf Acknowledgment.} The authors wish to thank A. Cossidente, A. Garcia,
J.W.P. Hirschfeld and G. Korchm\'aros for useful comments. The paper was
partially written while the first author was visiting IMECC-UNICAMP
supported by Capes-Brazil and Cnpq-Brazil.

\end{document}